\documentclass[submission]{dmtcs}

\usepackage[latin1]{inputenc}
\usepackage{subfigure}
\usepackage{mathtools}

\newcommand\polya{P\' olya}
\newcommand\Bin{{\rm Bin}}
\newcommand\Ber{{\rm Ber}}
\newcommand\Flajoletop{{\mathcal{M}}}
\newcommand\normal{{\mathcal{N}}}
\newcommand\convD{{\buildrel {\mathcal D} \over \longrightarrow}}
\newcommand\Prob{{\mathbb{P}}}

\newtheorem{prop}{Proposition}[section]
\newtheorem{theo}{Theorem}[section]
\newtheorem{example}{Example}[section]
\newtheorem{defin}{Definition}[section]
\newtheorem{rem}{Remark}[section]

\usepackage[round]{natbib}

\author{Basile Morcrette\addressmark{1}\thanks{Email: \email{Basile.Morcrette@inria.fr}. \\ Supported by ANR MAGNUM (10BLAN0204) and ANR BOOLE (09BLAN0011)}
  \and  Hosam M.~Mahmoud\addressmark{2}\thanks{Email: \email{hosam@gwu.edu}}}

\title[Analytic solvable balanced tenable urns with random entries]{Exactly solvable balanced tenable urns with random entries via the analytic methodology}
\address{\addressmark{1}Algorithms Project, INRIA Paris-Rocquencourt  \& LIP6, University Paris 6  (France) \\
  \addressmark{2} Department of Statistics, The George Washington University (U.S.A.)}

\received{18$^{th}$ January 2012}
\revised{\today}
\accepted{}

\begin{document}
\maketitle

\begin{abstract}
This paper develops an analytic theory for the study of some P\'olya urns with random rules. 
The idea is to extend the isomorphism theorem in Flajolet et al.~(2006), which connects deterministic balanced urns to a differential system for the generating function. 
The methodology is based upon adaptation of operators and  use of a weighted probability generating function.
Systems of differential equations are developed, and when they can be solved, they lead to characterization of the exact distributions underlying the urn evolution. We give a few illustrative examples.
\end{abstract}

\section{Introduction}
\label{sect:intro}

A classical definition of the \polya\ urn scheme specifies it as an urn which contains balls of up to two colors (say black and white), and which is governed by a set of evolution rules. A step then consists in randomly picking a ball from the urn, placing it back, and depending on its color, adding or discarding a fixed number of black and white balls.

\cite{FlDuPu06} introduces an analytic method to deal with a class of two-color \polya\ 
urns with deterministic addition rules. 
The method leads to a fundamental connection to differential equations via an isomorphism theorem. In the present manuscript, 
we extend the method to deal with \polya\ urns with \emph{random addition rules}. 
One can think of schemes with random replacements as mixtures of deterministic schemes, with appropriate probabilities.

We develop an analytic theory for the study of 
balanced tenable \polya\ urns with random rules (we shall make the terms ``balanced" and ``tenable" more precise in Section~\ref{sect:def}).
Balance and tenability are restrictions on the rules that admit an analytic treatment through generating functions.
Our main result is the following: {\emph{to each system of random rules satisfying the balance and tenability conditions, one can associate a differential system which describes 
the generating function counting the sequence of configurations of the urn scheme.}}
A formal statement of this theorem appears in Section~\ref{sect:theorem}.

When the system can be solved, it leads to characterization of the exact distributions underlying the urn evolution.
In Section~\ref{sect:examples}, we give some examples where the generating function is explicit. Some of them include uniform and binomial discrete random variable in their rules. Other examples investigate some coupon collector variants.
Finding the generating function and deducing the probability distributions is quite mechanical, owing to classical analytic combinatorics.

\section{Basic definitions}
\label{sect:def}

In \polya's classical urn model, we have an urn containing balls of up to two
different colors, say black and white. The system evolves with
regards to particular evolution rules: at each step, add (possibly a negative number) black and  white balls. 
These rules are specified by a $2\times 2$
matrix:
\begin{equation}
 \begin{pmatrix}
	\alpha     &  \beta \\
	\gamma &   \delta
\end{pmatrix},
  \qquad\alpha, \delta \in \mathbb{Z}\, , \quad \beta, \gamma\in\mathbb{Z}_{\geq 0}\, .
\end{equation}
The rows of the matrix are indexed with the colors black and white respectively.
The columns of the matrix are indexed by the same colors: they are from left to right indexed by black and white.

The asymptotics of this construct have been approached by traditional probabilistic methods in \cite{AtKa68}, \cite{Smythe96}
and \cite{Janson04a}. \cite{Mahmoud08} presents some of these findings.

In this manuscript we deal with a model where the constants 
$\alpha$, $\beta$, $\gamma$, 
and $\delta$ are
replaced by random variables.
We attempt to
provide an approach for the study of the urn's composition after
a finite, and possibly small, number of ball draws. This is an
equally important line of attack and can be viewed as more
important in practice,\footnote{For example, the Ehrenfest urn is a model for the exchange of gases. If the gas is being exchanged between an airconditioned room and the outside, the user is more interested in what happens within the next hour or so, not what ultimately happens at infinity.} when the number of draws is not sufficiently
large to warrant approximation by asymptotics.  We focus on exact distributions in this paper.

\subsection{Balanced tenable urns with random entries}
We consider urns for which the dynamics of change 
are now represented by
the replacement matrix
\begin{equation}
\begin{pmatrix}
	\mathcal{A}     &  \mathcal{B} \\
	\mathcal{C}     &  \mathcal{D}
\end{pmatrix},
\label{eq:urnrv}
\end{equation}
where $\mathcal{A}, \mathcal{B}, \mathcal{C}$ and $\mathcal{D}$ are discrete random
variables. If a random variable has a negative realization, it means we discard balls.    

Let $B_n$ and $W_n$ be respectively the number of black and white
balls after~$n$ draws from the urn. We call the pair $(B_n, W_n)$ the
\emph{configuration} of the urn after $n$ steps (draws).
We start with a deterministic initial configuration
$(B_0, W_0) = (\cramped{b_0}, \cramped{w_0})$. 
 At each step, a ball is uniformly drawn from the urn. We look at its color
and \emph{we put it back in the urn}: if the color is black, we add
$\mathcal{A}$ black balls and $\mathcal{B}$ white balls; if the color is white, we add $\mathcal{C}$
black and $\mathcal{D}$ white balls; the random variables are generated independently at each step.

We now recall a few definitions, that serve to delineate the scope of our study.

\begin{defin}[Balance]
  An urn with the replacement matrix (\ref{eq:urnrv}) is said to be \mbox{\emph{balanced}}, if the sums across
   rows are equal, that is if $\mathcal{A}+\mathcal{B} = \mathcal{C} + \mathcal{D}$, 
  and these sums are always constant and equal to a fixed value $\theta$. 
  The parameter~$\theta$ is called the \emph{balance} of the urn.
\end{defin}

All urns considered in this paper are balanced. As a consequence, the
random variables are governed by the relations
$\mathcal{B} = \theta - \mathcal{A}$, and $\mathcal{C} = \theta - \mathcal{D}$.
In addition, we will only
consider urns with nonnegative balance, $\theta \ge 0$. 
The case $\theta < 0$ refers to diminishing urn models. 
Questions and results on these models are completely different; 
see \cite{Hwang07} for an analytic development.

Since we allow the random variables $\mathcal{A}$ and $\mathcal{D}$ to take negative integer
values, it is necessary to consider the concept of tenability.

\begin{defin}[Tenability]
  An urn scheme is said to be \emph{tenable}, if it is
  always possible to apply a rule, i.e., if it never reaches
  a deadlocked configuration (that is, choosing with positive probability a rule that cannot be applied because there are not enough balls).
\end{defin}
To illustrate this notion, suppose we have an urn of black and
white balls, with the replacement matrix
\begin{displaymath}\begin{pmatrix}
	-2    &  2 \\
	 2 &   -2
\end{pmatrix}.\end{displaymath}
This matrix, together with an even initial number of black balls
and an even initial number of white balls form a tenable scheme.
Any other initial condition makes this scheme untenable. 
For instance, if we start with two black balls and one white ball, the scheme gets stuck
upon drawing a white ball.

\begin{rem}  
A tenable balanced urn scheme cannot have negative realizations for the nondiagonal entries 
$\mathcal{B} = \theta -\mathcal{A}$ and $\mathcal{C} = \theta-\mathcal{D}$, and they cannot have unbounded supports, either.
Indeed, if $\mathcal{B}$ has a negative realization, say $-k$ for positive $k$, we discard $k$ white balls upon withdrawing a black ball. 
In view of the balance condition, we must add at the same step $\theta +k$ black balls; we can apply this rule repeatedly until we take out all white balls, and the scheme 
comes to a halt (contradicting tenability).
Therefore, $\mathcal{B}$ and $\mathcal{C}$ must have only nonnegative realizations.
Furthermore, suppose the urn at some step has $N$ black balls; 
if $\mathcal{B}$ has an unbounded support, it includes a value greater than 
$\theta + N$. Here $\mathcal{A}$ must realize a value less then $-N$.  
A draw of a black ball requires taking out more than $N$ black
balls (contradicting tenability). \end{rem}

\begin{example}\label{ex:PolFried}
This example illustrates a mixed model urn in the class we are considering.
It has random entries in the replacement matrix and is balanced and tenable.
Consider the \polya-Friedman urn scheme with the 
replacement matrix 
\begin{equation}
\label{Eq:alternating}
\begin{pmatrix}
	\mathcal{B}_p     &  1-\mathcal{B}_p \\
	1-\mathcal{B}_p &   \mathcal{B}_p
\end{pmatrix},
\end{equation}
where $\mathcal{B}_p$ is a \Ber($p$) random variable.\footnote{The notation \Ber($p$) stands for the Bernoulli random variable with parameter $p$, that assumes the value $1$ with success probability $p$, and the value $0$ with failure probability $1-p$.} 
Such a scheme has balance one, and alternates between
\begin{equation*}
\begin{pmatrix}
1&0\\
0& 1\end{pmatrix}, \qquad \textmd{and}\qquad
\begin{pmatrix}
0&1\\
1& 0 \end{pmatrix},
\end{equation*}
that is between 
\emph{\polya-Eggenberger's urn}  (with probability $p$) and \emph{Friedman's urn}  (with probability 
$1-p$).
The transition in the limit distribution of black balls when $p$ goes from 0 to 1 is quite unusual. For $p=0$ (Friedman's urn), a properly normalized number
of black balls has the classical \emph{Gaussian distribution}, and for $p=1$  (\polya-Eggenberger's urn), a properly scaled number
of black balls has  a \emph{beta distribution}.
The scheme with replacement matrix~(\ref{Eq:alternating})
is an interesting example that warrants further investigation, and
will be a special case of two illustrative examples (Subsections~\ref{ssect:binom} and \ref{ssect:uniform}). 
\end{example}
\subsection{Definitions for the analytic methodology}
The first use of analytic combinatorics for the treatment of urn models is in \cite{FlGaPe05}, where some balanced urns are linked to a partial differential equation and elliptic functions.
We now describe the analytic tools we use for our study. For a more complete (yet accessible) account see \cite{Morcrette2012} (\S 2). For a precise overview of analytic combinatorics, see the reference book \cite{FlSe09}.

The total number of balls after $n$ draws, denoted by $\cramped{s_n:= B_n + W_n}$, is deterministic because of the balance condition. Indeed, at each step, we add a constant number of balls (the balance $\theta$), so 
\begin{equation}\label{total_number}
s_n = s_0 + \theta n = b_0 + w_0 + \theta n \, .
\end{equation}
The main tool in this study is the \emph{weighted probability
generating function}
\begin{equation}
\label{eq:pgf}
Q(x, y, z) := \sum_{n=0}^\infty \sum_{ b, w \ge 0}
        s_0 s_1\ldots s_{n-1}\, \Prob(B_n = b, W_n = w) \, x^b y^w \, \frac {z^n} {n!}\, .
\end{equation}
The variable $x$ counts black balls, $y$ counts white balls, and $z$ counts the number of draws. As $b$, $w$ and $n$ are governed by (\ref{total_number}), one of the variables is redundant. So we can set $y=1$ without loss of generality: it is sufficient to study a marginal distribution because it determines the joint probability distribution.

\begin{rem}
The kernel  $\cramped{s_0 s_1\ldots s_{n-1}}$ corresponds to the total number of paths of length $n$ starting from the initial configuration. Indeed, at step 0 we can pick a ball among the $s_0$ balls in the urn; at step 1 we can pick a ball among the $s_1$ balls in the urn, etc.
\end{rem}

\begin{rem} 
The coefficient\footnote{We use the operator $\cramped{[x_1^{j_1} \ldots x_k^{j_k}]}$ to extract the coefficient 
of $\cramped{x_1^{j_1} \ldots x_k^{j_k}}$ from a generating function of $\cramped{x_1, \ldots, x_k}$.
}  $\cramped{[x^b y^w z^n / n!]}Q(x,y,z)$ does not represent 
a combinatorial count; this quantity can be fractional. 
Nevertheless, it is closely connected to the notion 
of ``history" described in \cite{FlDuPu06}. A history is nothing but 
a finite sequence of draws from the urn. One history of length~$n$ describes one 
possible evolution from step $0$ to step $n$. 
It can be viewed as a path in the evolution tree of the urn. 
In the context of deterministic rules, 
histories are combinatorial counts. 
However, in our case, we need more information to have 
a complete description of a sequence of $n$ draws. 
At each step, we can choose different rules 
depending on the realization of the random entries at that step. 
So, we are now dealing with \emph{weighted histories}, 
that is the sequence of $n$ draws with a 
weight factor corresponding to the probability of applying the $n$ different rules 
we use along the path from step $0$ to step $n$. If we denote 
by $\cramped{Q_{n,b,w}}$ the contribution of weighted histories beginning 
at step $0$ in the configuration $(b_0,w_0)$ and ending at step $n$ in the \emph{configuration} $(b,w)$, we have
\begin{equation}\label{eq:comb_gf}
Q(x,y,z) = \sum_{n=0}^{\infty} \sum_{b, w \ge 0}  Q_{n,b,w}  \, x^b y^w \frac{z^n}{n!} \, .
\end{equation}
The kernel $\cramped{s_0 s_1 \ldots s_{n-1}}$ is obtained by adding all weighted histories of length $n$, that is by setting $x$ and $y$ to 1:
\begin{displaymath}
Q(1,1,z) = (1 - \theta z)^{- s_0 / \theta}  \, ;  \qquad  n!\, [z^n] Q(1,1,z) = s_0 s_1 \ldots s_{n-1} 
       = \theta^n \, \frac{\Gamma (n+s_0/\theta)}{\Gamma(s_0/\theta)} \, .
\end{displaymath}
In this way, we preserve the equivalence between the probabilistic model and the combinatorial 
aspects of weighted history. The balance condition is key to this equivalence:
\begin{equation}\label{eq:prob}
\Prob(B_n = b, W_n = w) = \frac{ \left[ x^b  y^w z^n  \right]  Q(x,y,z)}{\left[ z^n \right] Q(1,1,z)} = \frac{Q_{n,b,w}}{s_0 \ldots s_{n-1}}\, . 
\end{equation}
\end{rem}

\section{An isomorphism theorem for urn schemes with random entries}
\label{sect:theorem}
In this section, we state our main result through the following theorem which links the behavior of a tenable balanced urn to a differential system.
The theorem always yields a differential system, and the difficulty then lies in extracting information on the generating function: indeed, the system obtained may not be always explicitly solvable by currently known techniques.

\begin{theo}
\label{Theo:Flajolet}
Given a balanced tenable urn with
the replacement matrix
\begin{displaymath} \begin{pmatrix} 
\mathcal{A}   & \theta - \mathcal{A} \\
\theta - \mathcal{D}   &\mathcal{D} \end{pmatrix},\end{displaymath}
where $\mathcal{A}$ and $\mathcal{D}$ are discrete random variables with distributions
\begin{displaymath}
\Prob(\mathcal{A} = k) = \pi_k \, , \qquad \Prob(\mathcal{D} = k) = \tau_k\, ,
\end{displaymath}
for $-K \leq k \leq  \theta$, (for some $K>0$)
, the probability generating function of the urn is given by
\begin{displaymath}
Q\left(x,y,z\right) = X^{b_0}(z) \, Y^{w_0}(z) \, ,
\end{displaymath}
where $(b_0, w_0)$ is the starting configuration, and the pair $(X(t)\, , Y(t))$ is the solution to 
the differential system
\begin{displaymath} \left\{\begin{array}{ccc}
x'(t)  &=&   \displaystyle \sum_{k=-K}^\theta \pi_k \, x^{k+1}(t) \, y^{\theta-k}(t) \, ; \\
y'(t)  &=&    \displaystyle \sum_{k=-K}^\theta \tau_k \, x^{\theta-k}(t) \, y^{k+1}(t) \, ,
\end{array}\right.\end{displaymath}
applied at $t = 0$, 
and any initial conditions $x := x(0) \textrm{ and } y := y(0)$, such as
$x y \not = 0$.
\end{theo}


\smallskip

\proof
The proof follows and generalizes that in \cite{FlDuPu06}.
The basic idea in this proof 
is that the monomial $\cramped{x^i y^j}$ represents
the configuration $\cramped{(B_n, W_n) = (i,j)}$. 
The actions of the urn 
transform such a configuration into $\cramped{x^{i+k} y^{j+ \theta-k}}$, if a black ball is drawn and if the~$n$th realization of the
replacement matrix is such that $\mathcal{A} = k$, and $\theta - \mathcal{A} = \theta - k$ (which occurs with 
probability $\pi_k$); there are $i$ choices
for such a black ball.
Alternatively, the monomial is transformed 
into $\cramped{x^{i+\theta-k} y^{j+k}}$, if a white ball is drawn and if the $n$th realization
of the replacement matrix is such that $\theta - \mathcal{D} = \theta-k$, and $\mathcal{D} = k$ (which occurs with 
probability $\tau_k$); there are $j$ choices
for such a white ball. We represent the transformation by the operator
\begin{displaymath}
\Flajoletop := \sum_{k = - K}^{\theta} \pi_k \, x^k \, y^{\theta-k} \,  \Theta_x   + \tau_k  \, x^{\theta-k} \, y^k \,  \Theta_y \, ,
\end{displaymath}
where $\Theta_u  = u \partial_u$ is the \emph{pick and replace} operator.\footnote{See Subsection I.6.2 p.86 in \cite{FlSe09} for more details on the \emph{Symbolic method}.}
When this operator is applied to the configuration~$(i,j)$, it yields 
\begin{displaymath}
\Flajoletop(x^i y^j) = \sum_{k=-K}^\theta i \, \pi_k \, x^{i+k}\,  y^{j+\theta-k} + j \, \tau_k \, x^{i+\theta-k}\, y^{j+k} \, .
\end{displaymath}
This operator is used to represent all possible configurations of the urn at step $n+1$, 
given the configurations at step $n$. 
Indeed, if we write $Q(x,y,z) := \sum_{n=0}^\infty q_n (x,y) z^n / n!\, $, as a function in one variable~$z$, the operator $\Flajoletop$ represents the transition
\begin{displaymath} 
q_{n+1} (x,y) = \Flajoletop \left( q_n(x,y) \right) \quad \textmd{thus} \quad q_n (x,y) = \Flajoletop^n \left( q_0 (x,y)\right)  = \Flajoletop^n \left( \cramped{x^{b_0} y^{w_0}}\right).
\end{displaymath}
Let $x(t)$ and $y(t)$ be two functions that have Taylor series
expansion near $t=0$. Recall that the Taylor series expansion of the product 
$\cramped{x^\ell (t+z) \, y^m (t+z)}$ is
\begin{displaymath} 
x^{\ell} (t+z) \, y^{m}(t+z) = \sum_{n=0}^\infty  \frac {\partial^n} 
      {\partial t^n}  \left(x^{\ell} (t) \, y^{m}(t)\right) 
   \, \frac {z^n} {n!}.
\end{displaymath}
We also have
\begin{displaymath}
Q \left(x(t), y(t), z\right) = \sum_{n=0}^\infty \Flajoletop ^n \left( x^{b_0}(t)\, y^{w_0}(t)\right) \, \frac{z^n} {n!}.
\end{displaymath}
Take $\ell = b_0$, and $m = w_0$, and
these two expansions coincide, if the operators~$\Flajoletop$ 
and  $ {\partial} / {\partial t}$ are the same. This is possible if, for all $i, j \geq 0$,
\begin{eqnarray*}
\sum_{k=-K}^\theta i\,  \pi_k \, x^{i+k} (t)\,  y^{j +  \theta -k}(t) + j\,  \tau_k \, x^{i+\theta-k}(t) \, y^{j+k} (t)
   &   = & i\, x^{i-1} (t)\, y^j(t) \, x'(t) + j \, x^i (t)\, y^{j-1}  \, y'(t)  \, ,
\end{eqnarray*}
and can happen by choosing
\begin{displaymath}
\left\{\begin{array}{ccc}
x'(t)  &=&                       \displaystyle \sum_{k=-K}^\theta \pi_k \, x^{k+1}(t) \, y^{\theta-k}(t)  \, ;  \\
y'(t)  &=&                       \displaystyle \sum_{k=-K}^\theta \tau_k \, x^{\theta-k}(t) \, y^{k+1}(t)\, .
\end{array}\right.
\end{displaymath}
Hence, if $\left(X(t), Y(t)\right)$ is a solution of this system with initial conditions $x:=X(0)$ and $y:=Y(0)$,
we have $Q \left(X(t), Y(t), z\right) = \cramped{X^{b_0} (t+z)} \, 
\cramped{Y^{w_0}(t+z)}$.\footnote{A classic theorem of Cauchy and Kovalevskaya guarantees
the existence of a solution, and the condition $x y \neq 0$ guarantees an analytic expansion near $t=0$; see, for example, \cite{Folland95} (Chap.1 Sect.D., p. 46--55).}  
The statement follows by letting $t=0$.~\qed

\section{Examples of exactly solvable urns with random entries}
\label{sect:examples}
We now consider a variety of examples 
for which Theorem~\ref{Theo:Flajolet} admits an exact solution,
and  an exact probability distribution is obtained. 
In Subsection~\ref{ssect:coupon} we treat a variation of the coupon-collector problem.
In Subsections~\ref{ssect:binom} and \ref{ssect:uniform}, we study urns with binomial random variables and uniform random variables. In Subsection~\ref{ssect:mcoupon}, we study a variation of the coupon-collector problem that is modeled with three colors, demonstrating
that the analytic methods can be extended to more colors.

By probabilistic techniques,
\cite{Janson04a} and \cite{Smythe96} provide broad asymptotic urn
theories covering some of the urns in the class considered here; our
focus is on the \emph{exact distributions}.
\subsection{{An urn for coupon collection with delay}}
\label{ssect:coupon}
Collecting $k$ coupons can be represented by a number of urn schemes. 
One standard
urn model considers uncollected coupons to be balls of one color (say black)
and collected coupons to be balls of another color (say white); initially
all $k$ coupons are not collected (all $k$ balls are black). 
When a coupon is collected for the first time 
(a black ball is drawn from the urn),
one now considers that coupon type as acquired, 
so the black ball is recolored white
and deposited back in the urn. Collecting an already collected coupon type (drawing
a white ball from the urn) 
results in no change (the white ball is returned to the urn). Thus, the standard coupon
collection is represented by the replacement matrix
\begin{displaymath}
\begin{pmatrix}
-1  &  1\\
 0  &   0 \end{pmatrix}.
 \end{displaymath}
Occasionally, the coupon collector may misplace or lose a collected coupon \emph{immediately} after collecting it; the coupon still needs to be collected. That is,
the drawing of a black ball may sometimes result in no change too, delaying
the overall coupon collection.
This can be represented by the replacement matrix
\begin{equation}\label{eq:matrix_coupon}
\begin{pmatrix}
-\mathcal{B}_p   &\mathcal{B}_p\\
          0&   0 \end{pmatrix},
\end{equation}
where $\cramped{\mathcal{B}_p}$ is a \Ber($p$) random variable.
Under any nonempty starting conditions, this scheme is balanced and tenable.

\begin{prop}
\label{prop:coupon}
For the urn~(\ref{eq:matrix_coupon}), the exact probability distribution of the number $B_n$ of uncollected coupons (black balls)
after $n$ draws is given by
\begin{displaymath}
\Prob(B_n = b)   = \sum_{j=b}^{b_0} \, (-1)^{j-b} {\binom{b_0} {j}}  {\binom{j}{b}}   \left( \frac {s_0 -pj}{s_0}\right)^n .
\end{displaymath}
\end{prop}

\proof
The urn scheme~(\ref{eq:matrix_coupon}) is amenable to
the analytic method described, and leads to an exactly solvable system
of differential equations.
By Theorem~\ref{Theo:Flajolet}, the associated
differential system is
\begin{displaymath}\left\{\begin{array}{ccl}
 x'(t) &=& (1-p)\, x(t) +  p\, y(t)  \, ;\\ 
 y'(t) &=&  y(t) \, ,
\end{array}\right.\end{displaymath}
with solution
\begin{displaymath}
X(t)  =    y e^t + \left(x - y\right) e^{(1-p) t}\, , \qquad \textmd{and} \qquad
Y(t)  =    y e^t  \, .
\end{displaymath}

\noindent Thus,
\begin{displaymath}Q(x, 1, z) =  \bigl(e^z + (x-1) e^{(1-p) z}\bigr)^{b_0} e^{w_0 z}.\end{displaymath}
In the course of coupon collection the total number of coupons (balls in the urn)
does not change; that is, $\cramped{s_n = s_0}$, for all $n \geq 0$. Extracting coefficients, we get
\begin{eqnarray*}
\frac {s_0\ldots s_{n-1}} {n!} \, \Prob(B_n = b) &=& [x^b z^n] \, e^{w_0 z}  \sum_{j=0}^{b_0} { \binom{b_0} {j}}
              (x-1)^j e^{(1-p) j z} \times e^{(b_0 -j)z}  \\
        &=&  [z^n] \, e^{w_0 z}  \sum_{j=b}^{b_0} {\binom{b_0}{j}}(-1)^{j-b} {\binom{j}{b}}  e^{((1-p) j + b_0 -j)z} \, .
\end{eqnarray*}
Since $s_0 = b_0 + w_0$, we can write this as
\begin{eqnarray*}
\frac {s_0^n} {n!} \, \Prob(B_n = b)  &=& [z^n] \,  \sum_{j=b}^{b_0} (-1)^{j-b} {\binom{b_0}{j}} {\binom{j}{b}} e^{( s_0 -pj)z} \, .
\end{eqnarray*}
Consequently, we find the probability expression. \qed

\subsection{{A binomial urn}}
\label{ssect:binom}
Consider an urn scheme 
with replacement matrix
\begin{equation}
\label{eq:matrix_binom}
\begin{pmatrix} \mathcal{X}_{\theta,p}&\theta-\mathcal{X}_{\theta,p}\\
                  \theta-\mathcal{X}_{\theta,p}& \mathcal{X}_{\theta,p} \end{pmatrix},
\end{equation}
where $\cramped{\mathcal{X}_{\theta,p}}$ is distributed 
like $\Bin(\theta, p)$.\footnote{The 
notation $\Bin(k, p)$ stands for a binomially distributed random variable
that counts the number of successes in $k$ independent identically 
distributed trials with rate of success~$p$ per trial.}
Under any nonempty starting conditions, 
this is a balanced tenable scheme of balance $\theta$.
For an expression of the probability distribution of black balls after $n$ draws, 
it seems that  a solution for general $p$ is too difficult. 
We focus here on the tractable unbiased case $p= \frac 1 2$.
\begin{prop}
\label{prop:binom}
For the unbiased case $p=\frac1 2$ of the urn~(\ref{eq:matrix_binom}), the exact probability distribution of the number of black balls $B_n$ is 
\begin{displaymath}\Prob(B_n = b)  =  \frac{1} {2^{\theta n}} { \binom{\theta n}{b-b_0}} \, , \quad \textmd{\emph{for} } b \in \{b_0,\, b_0+1, \ldots,\, b_0 + \theta n\} .\end{displaymath}
\end{prop}

\noindent In other words, we have  $B_n = b_0 + \Bin \left(\theta n,  1/2 \right).$ 
Asymptotics follow easily from the normal approximation 
to the binomial distribution.
Let $\cramped{\normal(0, \sigma^2)}$ be a normal random variate with variance $\sigma^2$.
We see that
\begin{displaymath}
\frac{B_n - \frac 1 2 \theta n} {\sqrt n} \ \convD \ \normal\left( 0, \frac {\theta^2} {4}\right) \, ,
\end{displaymath}
where the symbol $\convD$ stands for convergence in distribution.

\begin{rem} 
The replacement matrix~(\ref{Eq:alternating}) 
is the special case where $\theta=1$, and as 
mentioned in the introduction
it alternates between \polya-Eggenberger's urn (with probability $p$)
and Friedman's urn (with probability $1-p$). So, one may think
of the urn process as a mixture model that chooses
between two given models at each step. It is well known that the number of
white balls in a pure \polya-Eggenberger urn (scaled by $n$) converges to a beta
distribution (see \cite{Polya31}), 
whereas the number of white balls in a pure Friedman's urn satisfies
the central limit tendency (see \cite{Freedman65}):
\begin{displaymath}
\frac{B_n - \frac 1 2 n} {\sqrt n} \ \convD \ \normal\left (0, \frac 1 {12}\right)\, .
\end{displaymath} 
In an unbiased mixture we get
$\normal (0,  1/4)$ as limit, 
with centering and scaling similar to that in Friedman's
urn. The normal limit for an unbiased mixture 
has a larger variance than a pure Friedman's urn, in view of occasional perturbation by the
\polya-Eggenberger choice. Because Friedman's urn is ``mean reverting," 
the entire mixed process is mean reverting, that is having tendency for average equilibrium around an even split.
We thus see that the Friedman effect is stronger than the \polya-Eggenberger.
\end{rem}

\proof[of Proposition~\ref{prop:binom}]
By Theorem~\ref{Theo:Flajolet}, the associated
differential system is
\begin{displaymath} \left\{
\begin{array}{ccc}
x'(t) &=& \displaystyle \sum_{k=0}^\theta  \, p^k (1-p)^{\theta-k}\, {\binom{\theta}{k}}\, x^{k+1}(t) \, y^{\theta-k}(t); \\ 
 y'(t) &=& \displaystyle \sum_{k=0}^\theta  \, p^k (1-p)^{\theta-k}\, {\binom{\theta}{k}}\, x^{\theta-k}(t) \, y^{k+1}(t).
\end{array}
\right.\end{displaymath}
The unbiased case, $p=  1/2$,  yields an explicit solution:
\begin{displaymath} 
X(t)  =    \frac {x} {\Bigl(1 - \theta \Bigl(\displaystyle \frac {x +y} 2 \Bigr)^\theta t \Bigr)^{1/\theta}} \, , 
\qquad \textmd{and}\qquad
Y(t)  =                     \frac {y} {\Bigl(1 - \theta \Bigl(\displaystyle \frac {x +y} 2 \Bigr)^\theta t \Bigr)^{1/\theta}} \, .
\end{displaymath}
The probability generating function $Q(x,y,z) = X(z)^{b_0}\ Y(z)^{w_0}$ gives
\begin{eqnarray*}
Q(x, y, z) &=&  \sum_{n=0}^\infty \, \sum_{ b, w \ge 0}
               s_0 s_1\ldots s_{n-1} \,  \Prob(B_n = b, W_n = w)\, x^b y^ w\, \frac {z^n} {n!}  \\
          &=&  \frac {x^{b_0} y^{w_0}}  {\Bigl(1 - \theta \Bigl(\displaystyle \frac {x + y} 2 \Bigr)^\theta z \Bigr)^{\frac {b_0+w_0} \theta}}.
 \end{eqnarray*}
 Note that $\Prob(B_n = b, W_n = w)$ is 0 for all values of $b$ and $w$, 
 except for $w = \theta n+ s_0 - b$,
 when the probability may differ from zero. In other words, the joint probability
 distribution can be determined from either marginal distribution. We compute
 \begin{eqnarray*}
Q(x, 1, z) &=&  \sum_{n=0}^\infty \, \sum_{ b, w \ge 0} \, 
               s_0 s_1\ldots s_{n-1}\, \Prob (B_n = b, W_n = w)\, x^b \, \frac {z^n} {n!}  \\
               &=& \sum_{n=0}^\infty  \, \sum_{b=0}^ \infty \,
               s_0 s_1\ldots s_{n-1}\, \Prob (B_n = b)\, x^b \, \frac {z^n} {n!}  \\
          &=&  \frac {x^{b_0}}  {\Bigl(1 - \theta \Bigl(\displaystyle \frac {x + 1} 2 \Bigr)^\theta z \Bigr)^{\frac {b_0+w_0} \theta}}.
 \end{eqnarray*} 
Extraction of coefficients yields\footnote{Recall that for $n,k \ge 0$ ,  ${\binom{-k}{n}}   =   \frac{-k (-k - 1) \ldots (-k + n -1)}{n!}$  .}
\begin{eqnarray*}
[x^b  z^n] \, Q(x,1, z) 
           &=& \frac {s_0\ldots s_{n-1}} {n!} \, \Prob (B_n = b) \\
           &=& [x^b  z^n] \,  x^{b_0}  \sum_{n=0}^\infty 
                      \Bigl(- \theta \Bigl(\displaystyle \frac {x + 1} 2 \Bigr)^\theta z \Bigr)^n     
                          {  \binom{{- \frac{b_0+w_0}{\theta}}} {n} }\\
             &=& [x^b] \,  x^{b_0}  \frac {1} {2^{\theta n}} \, (x + 1)^{ \theta n} \, \frac {s_0\ldots s_{n-1}} {n!}\\
             &=& [x^b] \,  x^{b_0}  \frac {1} {2^{\theta n}} \sum_{b=0}^{\theta n} {\binom{\theta n}{b}} \, x^b \, \frac {s_0\ldots s_{n-1}} {n!}.
\end{eqnarray*}
From this last expression we deduce the stated probability distribution for
the number of black balls.~~\qed

\subsection{{ A uniform urn}}
\label{ssect:uniform}
Consider an urn scheme 
with replacement matrix 
\begin{equation}
\label{eq:matrix_uniform}
\begin{pmatrix}
\mathcal{U}_\theta&\theta-\mathcal{U}_\theta \\
\theta-\mathcal{U}_\theta& \mathcal{U}_\theta \end{pmatrix},
\end{equation}
where $\mathcal{U}_\theta$ is a uniformly distributed random variable on the 
set $\{0, 1, \ldots, \theta\}$. Under any nonempty starting conditions, 
this is a balanced tenable scheme.

\begin{prop}
\label{prop:uniform}
The exact probability distribution of black balls for the urn~(\ref{eq:matrix_uniform}) is given by
\begin{displaymath}
\Prob(B_n = b)  =  \frac{T_{\theta, n, b-b_0}} {(\theta+1)^n} \, , \quad \textmd{\emph{for} } b \in \{b_0, \, b_0+1, \ldots, \,  b_0 + \theta n\} , 
\end{displaymath}
where $\cramped{T_{\theta, n, k}} := [x^k] \, (1+ x + \ldots + x^\theta)^n$.
\end{prop}

The numbers $\cramped{T_{\theta, n, k}}$
are known in the classical literature (see \cite{Euler01}) and have numerous combinatorial interpretations. 
For instance, in \cite{FlSe09} (Subsection I.15, page 45), $T_{r-1, k , n-k}$ is the number of compositions of size $n$ with $k$ summands each at most $r$.
In other words, $\cramped{T_{2,n,k}}$, called the trinomial coefficient, is the number of distinct ways in which $k$ indistinguishable
balls can be distributed over $n$ distinguishable urns allowing 
at most two balls to fall in each urn.
Several other interpretations related to counting strings,  
partitions and unlabeled trees of height 3 exist; see \cite{Andrews90}.

\proof[of Proposition~\ref{prop:uniform}]
By Theorem~\ref{Theo:Flajolet}, the associated
differential system is
\begin{displaymath}\left\{\begin{array}{ccc}
 x'(t) &=& \displaystyle \sum_{k=0}^\theta \, \frac 1 {\theta+1}\, x^{k+1}(t) \, y^{\theta-k}(t)  \, ;\\ 
 y'(t) &=& \displaystyle \sum_{k=0}^\theta  \, \frac 1 {\theta+1}\, x^{\theta-k}(t) \, y^{k+1}(t)  \, .
\end{array}\right.\end{displaymath}
This yields an explicit solution:
\begin{displaymath}
X(t)  =                    \frac {x} {\Bigl(1 - \displaystyle \frac \theta {\theta+1} \Bigl( \sum_{l=0}^\theta
          x^l y^{\theta-l}\Bigr)t \Bigr)^{1/\theta}}   \, ; \qquad
Y(t)  =                    \frac {y} {\Bigl(1 - \displaystyle \frac \theta {\theta+1} \Bigl( \sum_{l=0}^\theta
          x^l y^{\theta-l} \Bigr)t \Bigr)^{1/\theta}}   \, .
\end{displaymath}
Thus, the probability generating function $Q(x,1,z) = X(z)^{b_0}\, Y(z)^{w_0}$ is
 \begin{eqnarray*}
Q(x, 1, z)  &=& \sum_{n=0}^\infty \sum_{b=0}^ \infty
                       s_0 s_1\ldots s_{n-1}\, \Prob(B_n = b)\, x^b \, \frac {z^n} {n!}  \\
                &=&  \frac {x^{b_0} }  {\Bigl(1 - \displaystyle \frac \theta {\theta+1} \Bigl( \sum_{l=0}^\theta
                    x^l\Bigr) z \Bigr)^{\frac {b_0+w_0} \theta}}  \, .
 \end{eqnarray*} 
Extraction of coefficients yields
\begin{eqnarray*}
\frac {s_0\ldots s_{n-1}} {n!} \, \Prob(B_n = b) &
            = & [x^b  z^n] \,  x^{b_0}  \sum_{n=0}^\infty \Bigl(-  \displaystyle \frac \theta {\theta+1} \Bigl( \sum_{l=0}^\theta
                    x^l\Bigr) z \Bigr)^n     
                          {\binom{ - \frac {b_0+w_0} \theta}{n}} \\
        &  = &  [x^b] \,  x^{b_0}  \frac {1} {(\theta+1)^n}\left(1+ x + \cdots + x^\theta\right)^n \,
                      \frac {s_0\ldots s_{n-1}} {n!}.
\end{eqnarray*}
We thus obtain the stated probability expression.\qed

\subsection{{An urn for a two-type coupon collection}}
\label{ssect:mcoupon}
The method has obvious extensions for $k\geq 2$ colors. 
Consider, for instance, the balanced tenable three-color coupon 
collection urn scheme
with the entries
\begin{equation}
\begin{pmatrix}
 -1& \mathcal{B}_p& 1-\mathcal{B}_p\\
                  0&0&0  \\
                  0&0&0 \end{pmatrix},
                  \label{Eq:threecolor}
\end{equation}
in which a collected coupon is randomly categorized to fall in one of two classes.
Here again, $\cramped{\mathcal{B}_p}$ is a \Ber($p$) random variable.
The uncollected coupons correspond to black balls in the urn (replacement
rules on the first row), and the white balls of the previous example are now ramified
into red and green balls (second and third rows). The columns of the matrix are indexed
by the same colors: they are from left to right indexed by black, red and green.
Let $\cramped{B_n}$, $\cramped{R_n}$, $\cramped{G_n}$ be respectively 
the number of black, red and green balls
in the urn after $n$ draws.
 
 \begin{prop}
 \label{prop:mcoupon}
 For the urn~(\ref{Eq:threecolor}), the distribution of red balls $R_n$ is given by the following probabilities:
 \begin{displaymath}
 \Prob(R_n = r) = \left(\frac {p}{(1-p)}\right)^r \sum_{j=0}^{b_0} (-1)^j {\binom{j}{r}} {\binom{b_0}{j}} (1-p)^j\sum_{k=0}^j {\binom{j}{k}}
   (-1)^k  \left(\frac k {s_0}\right) ^n. 
   \end{displaymath}
 \end{prop}

\proof 
This scheme has the underlying differential system
\begin{displaymath}\left\{\begin{array}{ccl}
 x'(t) &=& p \, y(t) + (1-p)\, h(t)  \,  ;\\ 
 y'(t) &=& y(t)   \,  ;\\
 h'(t)  &=& h(t) \, .
\end{array}\right.\end{displaymath}
This system of differential equations has the solution
\begin{eqnarray*}
X(t) &=& \left(p\, y + (1-p)\, h\right) (e^t - 1) + x   \,  ; \\
Y(t) &=& y\, e^t   \,  ; \\
H(t) &=& h\, e^t   \, .
\end{eqnarray*}
Note that in this scheme the total number of balls after any number
of draws remains the same at all times.
The corresponding generating function is
\begin{displaymath}Q(x, y, h, z) = \sum_{n=0}^\infty \ \sum_{b, r, g \ge 0}
        s_0 ^n \, \Prob(B_n = b, R_n = r, G_n = g) x^b y^r h^g\, \frac {z^n} {n!}  \, .\end{displaymath}
The extended isomorphism theorem then gives
\begin{displaymath}Q\left(x, y, h, z\right) = X^{b_0}(t) \, Y^{r_0}(t)\, H^{g_0} (t)  \, .\end{displaymath}
One can extract univariate and joint distributions from this.
For instance,
\begin{eqnarray*}
\frac {s_0^n} {n!}\, \Prob(R_n = r) &=& \left[y^r z^n\right] \, Q(1, y, 1, z) \\
         &=& \left[y^r z^n\right] \, \left( p \, (y-1)\, \left(e^z - 1\right)+e^z \right)^{b_0}  y^{r_0} e^{r_0 z} e^{g_0 z}  \, .  
\end{eqnarray*}
Extracting coefficients (details omitted) we find the probability distribution for
red balls. \qed

\section{Conclusion}   
We presented a generalization of the isomorphism theorem for deterministic schemes.
The generalization covers the class of balanced tenable
two-color urn schemes with random entries.
The result is particularly useful when the system of ordinary differential equations
is amenable to a simple solution, from which exact probability distributions
can be extracted in a fairly mechanical fashion.

\begin{figure}[htbp]
\begin{center}
\subfigure[$p = 0$\label{grap0}]{\includegraphics[height=2.2cm]{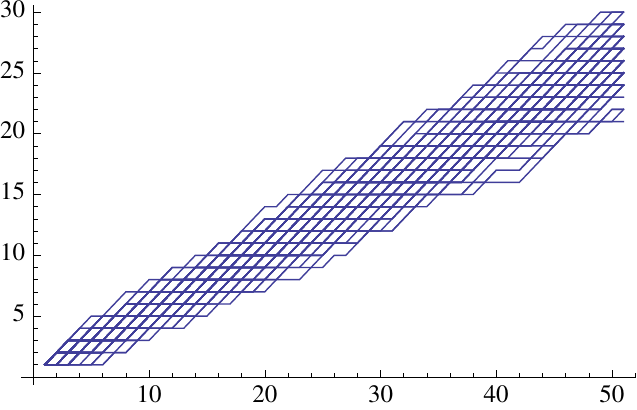}}
\subfigure[$p = 0.4$\label{grap04}]{\includegraphics[height=2.2cm]{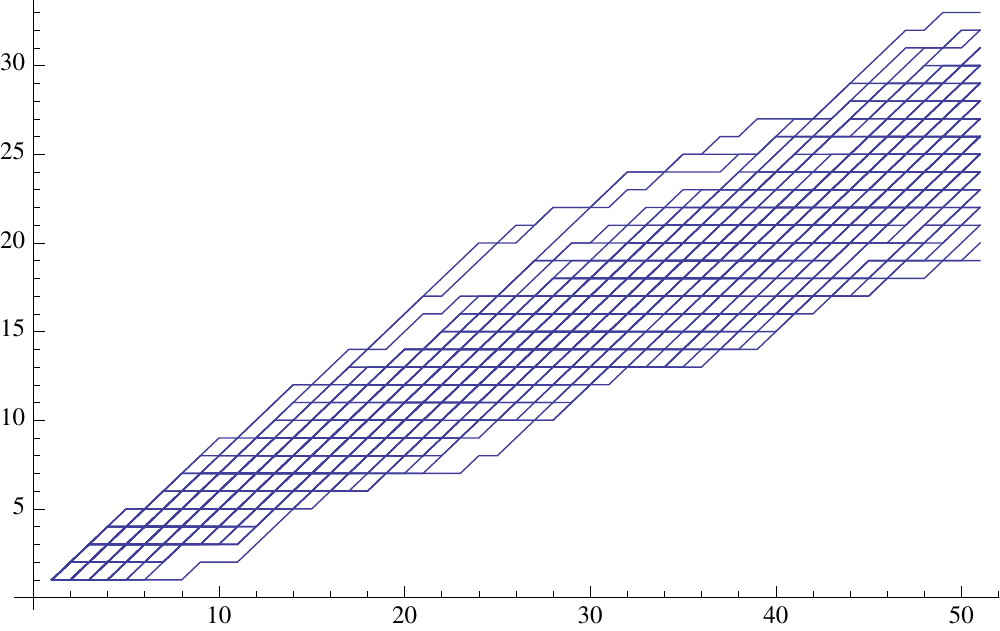}}
\subfigure[$p = 0.8$\label{grap08}]{\includegraphics[height=2.2cm]{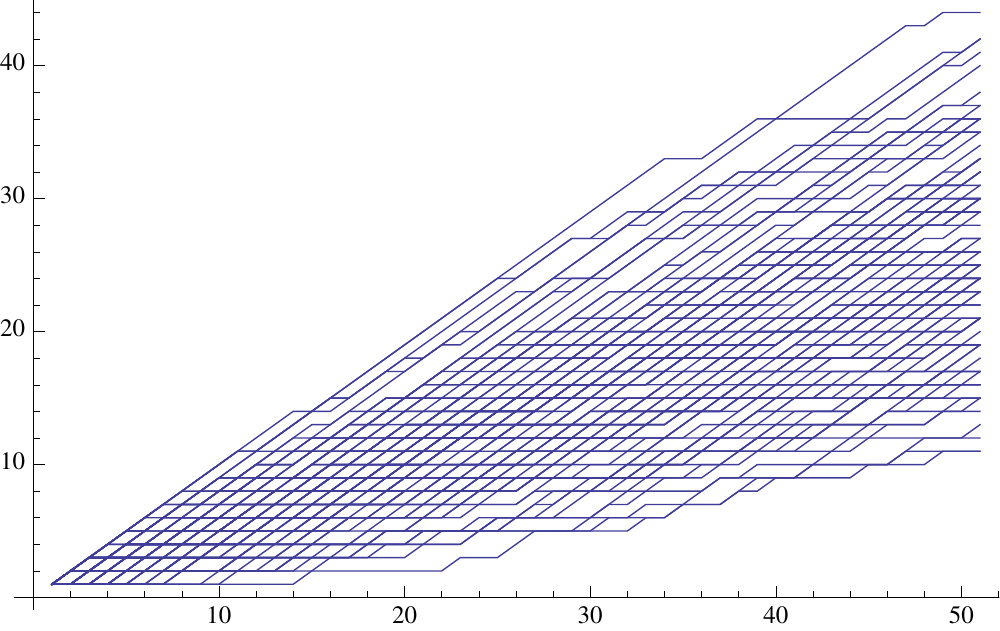}}
\subfigure[$p = 1$\label{grap1}]{\includegraphics[height=2.2cm]{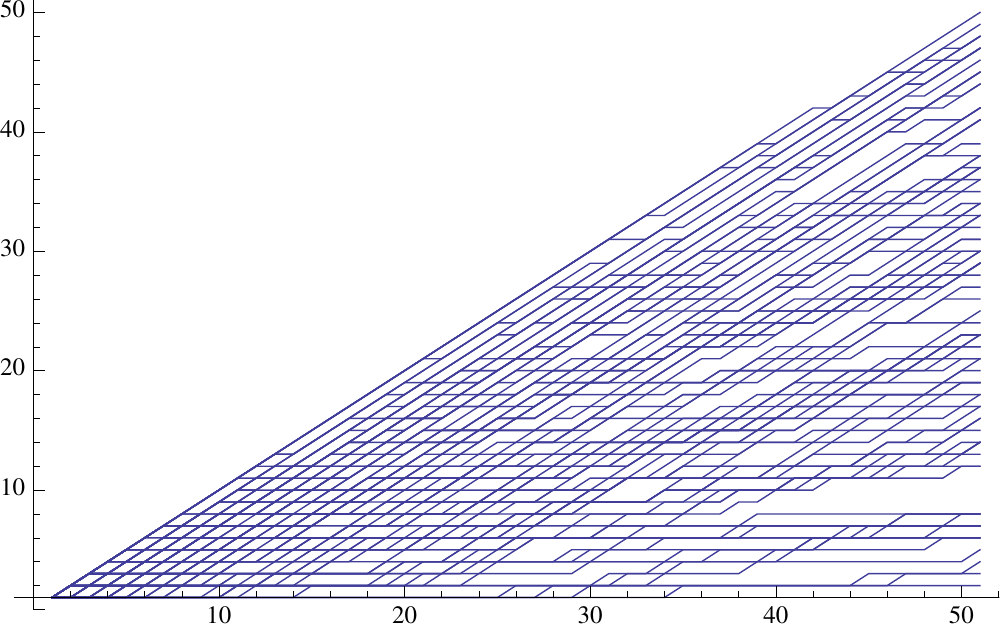}}
\caption{Each graph represents a simulation of 100 histories of length 50 for the \polya--Friedman's urn (\ref{Eq:alternating}) (see Example~\ref{ex:PolFried}) with parameter $p \in \{ 0,\, 0.4,\, 0.8,\, 1\} .$ There is a transition between a \emph{Gaussian}~behavior (Fig.\ref{grap0}) and a \emph{beta}~behavior (Fig.\ref{grap1}).}
\label{PolFried}
\end{center}
\end{figure}

\begin{figure}[htbp]
\begin{center}
\includegraphics[height=5.1cm]{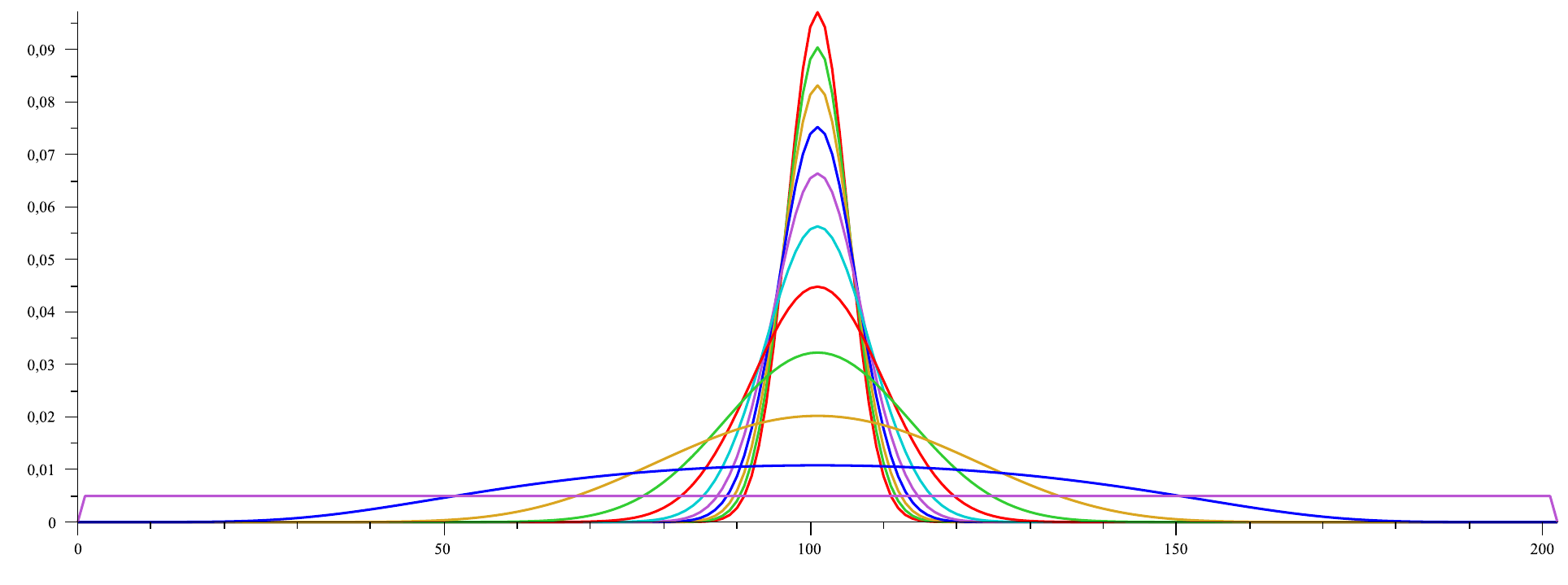}

\caption{
Normalized distribution of the number of black balls after 200 draws ($B_{200}$) for the \polya--Friedman's urn~(\ref{Eq:alternating}) with different values of $p \in \{ 0,\, 0.1,\, 0.2,\, 0.3,\, 0.4,\, 0.5,\, 0.6,\, 0.7,\, 0.8,\, 0.9,\, 1 \}$.
Here, the starting configuration is one black and one white ball. We observe a transition between a Gaussian distribution and a uniform distribution (which is a special case of beta distribution, when $b_0=w_0=1$).}
\label{fig:gauss}
\end{center}
\end{figure}

Our last example (see Proposition~\ref{prop:mcoupon}) 
shows it is straightforward to extend our theorem to additional colors: an urn with $k$ colors will lead to a differential system of $k$~equations. 
This theorem invites further work, using asymptotic tools from \cite{FlSe09}, to obtain limit laws.
For instance, the first example (see Proposition~\ref{prop:binom}, a model mixing \polya\  and Friedman urn models, and Fig.~\ref{PolFried} and~\ref{fig:gauss})  raises various questions such as: Can we obtain an exact explicit expression for the generating function for general $p \in (0,1)$? What is the limit law, and what is the asymptotic behavior? 
How does the phase transition between a Gaussian regime and a beta regime take place?


\bibliography{dmtcs-aofa12}

\begin{thebibliography}{14}
\providecommand{\natexlab}[1]{#1}
\providecommand{\url}[1]{\texttt{#1}}
\expandafter\ifx\csname urlstyle\endcsname\relax
  \providecommand{\doi}[1]{doi: #1}\else
  \providecommand{\doi}{doi: \begingroup \urlstyle{rm}\Url}\fi

\bibitem[Andrews(1990)]{Andrews90}
G.~Andrews.
\newblock Euler's ``exemplum memorabile inductionis fallacis" and
  $q$--trinomial coefficients.
\newblock \emph{Journal of the American Mathematical Society}, 3:\penalty0
  653--669, 1990.

\bibitem[Athreya and Karlin(1968)]{AtKa68}
K.~B. Athreya and S.~Karlin.
\newblock Embedding of urn schemes into continuous time {M}arkov branching
  processes and related limit theorems.
\newblock \emph{Annals of Mathematical Statistics}, 39:\penalty0 1801--1817,
  1968.

\bibitem[Euler(1801)]{Euler01}
L.~Euler.
\newblock De evolutione potestatis polynomialis cuiuscunque $(1+x+x^2+ x^3+x^4+
  \mbox{\ldots})^n$.
\newblock \emph{Nova Acta Academiae Scientarum Imperialis Petropolitinae},
  12:\penalty0 47--57, 1801.

\bibitem[Flajolet and Sedgewick(2009)]{FlSe09}
P.~Flajolet and R.~Sedgewick.
\newblock \emph{Analytic Combinatorics}.
\newblock Cambridge University Press, 2009.

\bibitem[Flajolet et~al.(2005)Flajolet, Gabarr\'o, and Pekari]{FlGaPe05}
P.~Flajolet, J.~Gabarr\'o, and H.~Pekari.
\newblock Analytic urns.
\newblock \emph{Annals of Probability}, 33:\penalty0 1200--1233, 2005.

\bibitem[Flajolet et~al.(2006)Flajolet, Dumas, and Puyhaubert]{FlDuPu06}
P.~Flajolet, P.~Dumas, and V.~Puyhaubert.
\newblock Some exactly solvable models of urn process theory.
\newblock In P.~Chassaing, editor, \emph{Fourth Colloquium on Mathematics and
  Computer Science}, volume~AG of \emph{DMTCS Proceedings}, pages 59--118,
  2006.

\bibitem[Folland(1995)]{Folland95}
G.~Folland.
\newblock \emph{Introduction to Partial Differential Equations}.
\newblock University Press, Princeton, New Jersey, 1995.

\bibitem[Freedman(1965)]{Freedman65}
D.~A. Freedman.
\newblock Bernard {F}riedman's urn.
\newblock \emph{Annals of Mathematical Statistics}, 36:\penalty0 956--970,
  1965.

\bibitem[Hwang et~al.(2007)Hwang, Kuba, and Panholzer]{Hwang07}
H.-K. Hwang, M.~Kuba, and A.~Panholzer.
\newblock Analysis of some exactly solvable diminishing urn models.
\newblock In \emph{Formal Power Series and Algebraic Combinatorics (FPSAC),
  Tianjin, China}, 2007.
\newblock 12 pages.

\bibitem[Janson(2004)]{Janson04a}
S.~Janson.
\newblock Functional limit theorems for multitype branching processes and
  generalized {P}\'olya urns.
\newblock \emph{Stochastic Processes and Applications}, 110\penalty0
  (2):\penalty0 177--245, 2004.

\bibitem[Mahmoud(2008)]{Mahmoud08}
H.~M. Mahmoud.
\newblock \emph{Urn Models}.
\newblock Chapman, Orlando, Florida, 2008.

\bibitem[Morcrette(2012)]{Morcrette2012}
B.~Morcrette.
\newblock Fully analyzing an algebraic {P}\'olya urn model.
\newblock In \emph{LATIN 2012}, volume 7256 of \emph{Lecture Notes in Computer
  Science}, pages 568--581, 2012.

\bibitem[P\'olya(1931)]{Polya31}
G.~P\'olya.
\newblock Sur quelques points de la th\'eorie des probabilit\'es.
\newblock \emph{Annales de l'Institut Poincar\'e}, 1:\penalty0 118--161, 1931.

\bibitem[Smythe(1996)]{Smythe96}
R.~T. Smythe.
\newblock Central limit theorems for urn models.
\newblock \emph{Stochastic Processes and their Applications}, 65\penalty0
  (1):\penalty0 115--137, 1996.

\end{thebibliography}
\bibliographystyle{abbrvnat}

\end{document}